\documentclass[12pt,a4paper,twoside,reqno]{amsart}

\usepackage{amsmath}
\usepackage{amssymb}
\usepackage{amsfonts}
\usepackage{amsthm}
\usepackage{amscd}
\usepackage{latexsym}
\usepackage[margin=3cm,heightrounded=true,centering]{geometry}
\usepackage[pdftex]{hyperref} 

\newcommand{\enveq}[1]{\begin{equation}#1\end{equation}}

\newcommand{\bma}{\left(\begin{array}{cc}}
\newcommand{\ema}{\end{array}\right)}
\newcommand{\bca}{\left(\begin{array}{c}}
\newcommand{\eca}{\end{array}\right)}

\newcommand{\cas}{C_{q}}

\numberwithin{equation}{section} 

\theoremstyle{plain} 
\newtheorem{thm}{Theorem}[section]
\newtheorem{lemma}[thm]{Lemma}
\newtheorem{prop}[thm]{Proposition}

\newtheorem{remark}[thm]{Remark}

\theoremstyle{definition} 
\newtheorem{defn}{Definition}[section]

\newcommand{\cd}{\cdot}
\newcommand{\clc}{\cdot\ldots\cdot}
\newcommand{\ot}{\otimes}
\newcommand{\op}{\oplus}

\newcommand{\olo}{\otimes\ldots\otimes}

\newcommand{\ci}{\circ}

\newcommand{\nn}{\mathbb{N}}
\newcommand{\zz}{\mathbb{Z}}

\newcommand{\cc}{\mathbb{C}}
\newcommand{\al}{\alpha}

\newcommand{\De}{\Delta}
\newcommand{\ga}{\gamma}

\newcommand{\la}{\lambda}

\newcommand{\pa}{\partial}
\newcommand{\si}{\sigma}
\newcommand{\Si}{\Sigma}
\newcommand{\te}{\theta}

\newcommand{\ze}{\zeta}
\newcommand{\C}[1]{\mathcal{#1}}

\newcommand{\T}[1]{\textrm{#1}}
\newcommand{\E}[1]{\emph{#1}}

\newcommand{\fork}[2]{\left\{ \begin{array}{#1} #2 \end{array} \right.}

\newcommand{\mat}[2]{\left(\begin{array}{#1} #2 \end{array} \right)}

\newcommand{\q}{\qquad}
\newcommand{\qq}{\qquad \qquad}

\newcommand{\su}{\subseteq}

\newcommand{\SU}{SU_{q}(2)}

\title{A twisted spectral triple for quantum $SU(2)$}

\author{Jens Kaad}
\address{Hausdorff Center for Mathematics,
Universit\"at Bonn,
Endenicher Allee 60,
53115 Bonn,
Germany}
\email{jenskaad@hotmail.com}

\author{Roger Senior}
\address{Mathematical Sciences Institute, Australian National University, Acton, ACT, 0200, Australia}
\email{roger.senior@anu.edu}

\thanks{The first author was supported by the Hausdorff Center for Mathematics (Bonn) and the Australian Research Council.}

\begin{document}

%

\maketitle

\vspace{15pt}

\begin{abstract}
We initiate the study of a $q$-deformed geometry for quantum $SU(2)$. In contrast with the usual properties of a spectral triple, we get that only twisted commutators between algebra elements and our Dirac operator are bounded. Furthermore, the resolvent only becomes compact when measured with respect to a trace on a semifinite von Neumann algebra which does not contain the quantum group. We show that the zeta function at the identity has a meromorphic continuation to the whole complex plane and that a large family of local Hochschild cocycles associated with our twisted spectral triple are twisted coboundaries.
\end{abstract}

\vspace{60pt}

\tableofcontents

\newpage

\section{Introduction}
In this paper we present an example of a $q$-deformed geometry for the quantum group version of the classical Lie group $SU(2)$. Our example differs from the spectral triples constructed by Chakraborty and Pal and by Dabrowsky \emph{et al}, \cite{chakpal,dablandsitsuijvar}, in a substantial way. First of all, the spectrum of our $q$-Dirac operator is not a linear function of $n \in \nn$, instead the dependency is exponential with slope given by the deformation parameter $\log(q)$. In particular we expect that invariants constructed on the basis of the $q$-Dirac operator will be sensitive to the $q$-deformation in an interesting way. On a more conceptual level, our geometrical framework is not captured by the existing notions of spectral triples. To be more precise, we get that the commutator between the $q$-deformed Dirac operator and algebra elements only becomes bounded when twisted by an automorphism. Furthermore, the resolvent is not compact when measured by the usual operator trace, instead the resolvent lies in the compacts of a semifinite von Neumann algebra. We will refer to our data as a twisted spectral triple.

The presence of twisted bounded commutators in the quantum group setting has been predicted by Connes and Moscovici in their paper \cite{conmosc}. Their main examples are obtained via conformal equivalence and foliations but they end the paper with the following comment on the quantum group case:

\begin{quote}
"The domain of quantum groups is a natural arena where twisting frequently occurs and where the [...] extension of the notion of spectral triple could be useful."
\end{quote}

Let us recall the Connes-Moscovici extension of a spectral triple.

\begin{defn}\label{sispectral}
Let $A$ be a unital $C^*$ algebra and let $\C A \su A$ be a dense $*$-subalgebra which comes equipped with an automorphism $\sigma \in \T{Aut}(\C A)$ such that $\sigma(x^*) = (\sigma^{-1}(x))^*$ for all $x \in \C A$. Then a \emph{$\sigma$-spectral triple} $(\C A, \C H, D)$ is given by a representation of $A$ on the Hilbert space $\C H$, while $D$ is a selfadjoint operator such that
\begin{enumerate}
\item The resolvent $(i + D)^{-1} \in \C K(\C H)$ lies in the compacts of the Hilbert space $\C H$.
\item The twisted commutator $[D, x]_{\si} := D x - \si(x) D$ extends to a bounded operator on $\C H$ for all $x \in \C A$.
\end{enumerate}
A graded $\si$-spectral triple is similarly defined, with the addition of a grading operator $\gamma \in \C L(\C H)$ satisfying $\gamma^* = \gamma$, $\gamma^{2} = I$, and which commutes with the action of $A$ but anti-commutes with $D$. A $\sigma$-spectral triple is said to be \emph{Lipschitz-regular} if the additional condition $[|D|, x]_{\sigma} \in \C L(\C H)$, for all $x \in \C A$ is satisfied.
\end{defn}

So far however, the main examples of spectral triples over $q$-deformed Lie groups are of a non-twisted nature, and the Dirac operators of these examples share a lot of properties with the classical Dirac operators. In particular the corresponding index invariants are unable to detect the $q$-deformation. See \cite{nestusI,nestusII}.

A Lipschitz-regular $\sigma$-spectral triple $(\C A, \C H, D)$ defines an ordinary $K$-homology class $[\C H, D(1 + D^2)^{-1/2}] \in K^{\bullet}(A)$, and therefore pairs with $K$-theory via the index pairing. Furthermore, under the additional assumption of finite summability the Chern-Connes character of the spectral triple defines a class in cyclic cohomology, and the pairing of this cocycle with the Chern character of $K$-theory classes computes the index pairing.

It is tempting to extend the notion of a $\sigma$-spectral triple to the semifinite case by replacing the Hilbert space compacts by the compacts with respect to some semifinite trace. Indeed, this idea has been succesful in the case of an ordinary spectral triple leading to new versions of the Connes-Moscovici local index theorem, see \cite{carphilrensukI,carphilrensukII,conmos}. We would like to emphasize that such an extension of the concept does not incorporate the example which we are looking at here in a sensible way. In order to construct a class in $KK$-theory from a semifinite spectral triple it is a necessary condition that the algebra is contained in the semifinite von Neumann algebra, see \cite{KNR}. In our example, the resolvent of the $q$-Dirac operator is compact only in the fixed point von Neumann algebra for a modular action that is non-trivial on the coordinate algebra for quantum $SU(2)$. This observation indicates that the notion of a $\si$-spectral triple should be combined with some of the modular ideas appearing in the recent paper by Rennie and the second author, see \cite{RS}. See also the paper \cite{carnesnesren} for an investigation of similar problems. We defer the proper treatment of these matters until more examples have been discovered.

Apart from proving the basic properties of our twisted spectral triple, we will also consider the zeta function associated with the $q$-Dirac operator and the semifinite trace. In particular, we show that this zeta function has a meromorphic extension to the whole complex plane. Finally, we investigate a rather large family of twisted local Hochschild cocycles arising from our data. These local Hochschild cocycles are constructed using twisted commutators with the $q$-Dirac operator and applying some twisted trace. We show that our local Hochschild cocycles are twisted coboundaries no matter which twisted trace we choose. This leaves the understanding of the homological dimension of our example open. And more generally, it leaves the search for a good twisted cohomological invariant of the $q$-deformed geometry unsettled.

\section*{Acknowledgements} The authors would like to thank Matthias Lesch, Ryszard Nest and Adam Rennie for many helpful discussions.

\section{Notation and conventions}\label{notcon}
In this section we fix some notation and conventions for the quantum group versions of the classical Lie-group $SU(2)$. These $q$-deformations of the special unitary matrices were introduced by Woronowicz in \cite{woron}, and we denote them by $SU_q(2)$.

Let us fix some real number $q \in (0,1)$. We start by introducing the coordinate algebra $\C A:= \C O(\SU)$ for the quantum group $SU_q(2)$. This is the unital $\cc$-algebra generated by the symbols $a,b,c,d \in \C A$ with relations
\[
\begin{split}
& ab = qba \, \, , \, \, ac = q ca 
\, \, , \, \,  bd = q db 
\, \, , \, \,  cd = q dc \, \, , \, \, bc = cb \\
& ad - q bc = 1 \, \, , \, \,  da - q^{-1} bc = 1.
\end{split}
\]
The $\cc$-algebra $\C A$ becomes a $*$-algebra when equipped with the involution defined by $* : \C A \to \C A$, $a^* = d$ and $b^* = -q c$. The coordinate algebra $\C A$ can be given the structure of a Hopf-$*$-algebra, see \cite{KS}. We will however not make use of this fact.
%
%

We let $\C U:= \C U_q(su_2)$ be the unital $\cc$-algebra generated by the symbols $e,f,k,k^{-1} \in \C U$ with relations
\begin{equation}\label{eq:relU}
\begin{split}
& k^{-1}k = 1 = kk^{-1} \\
& ke = q ek \, \, , \, \,  kf = q^{-1}fk \\
& ef - fe = (k^2 - k^{-2})/ (q - q^{-1})
\end{split}
\end{equation}
The $\cc$-algebra $\C U$ becomes a $*$-algebra when equipped with the involution given by $* : \C U \to \C U$, $e^* = f$ and $k^* = k$. The $*$-algebra $\C U$ can also be given the structure of a Hopf-$*$-algebra, see \cite{KS}.


We define the element
\begin{equation}\label{eq:casimir}
\begin{split}
c_q 
& := ef + (q^{-1} - q)^{-2} (q^{1/2} k^{-1} - q^{-1/2} k)^2 \\
& = fe + (q^{-1} - q)^{-2} (q^{-1/2} k^{-1} - q^{1/2} k)^2
\end{split}
\end{equation}
in $\C U$ and refer to it as the (quantum) Casimir. This element generates the centre of $\C U$.

There exists an algebra homomorphism $\pa : \C U \to \T{End}(\C A)$ which we will refer to as the \emph{left action}. For the precise definition we refer to \cite{KS}. We note that the action of the generator $k$ yields an automorphism $\pa_k \in \T{Aut}(\C A)$ whereas the actions of $e$ and $f$ satisfy the derivational rule
\begin{equation}\label{eq:twistdev}
\begin{split}
\pa_e(xy) & = \pa_e(x) \pa_{k^{-1}}(y) + \pa_k(x) \pa_e(y) \\
\pa_f(xy) & = \pa_f(x) \pa_{k^{-1}}(y) + \pa_k(x) \pa_f(y)
\end{split}
\end{equation}
for all $x,y \in \C A$. The automorphism $\si_L := \pa_{k^2} \in \T{Aut}(\C A)$ will play a special role in the present text and we will refer to it as the \emph{left modular automorphism}. It can be proved that the left modular automorphism satisfies the unitarity condition $\si_L(a^*) = \si_L^{-1}(a)^*$ of Definition \ref{sispectral}.
%
%
%

The quantum group $SU_q(2)$ comes equipped with a Haar-state which we will denote by $h : SU_q(2) \to \cc$. The Haar-state is a faithful state \cite{woron}. We will use the notation $\C H_h$ for the associated $GNS$-space which then comes equipped with an action $\pi : \SU \to \C L(\C H_h)$ of $SU_q(2)$ as left-multiplication operators. Furthermore, using the left action $\pa : \C U \to \T{End}(\C A)$ we get a representation of $\C U$ as unbounded operators on $\C H_h$ with domain given by the coordinate algebra $\C A \su \C H_h$. Under this representation, the $*$-operation corresponds to taking the \emph{formal} adjoint. For an element $g \in \C U$ we will use the capital letter $G : \C A \to \C H_h$ for the associated unbounded operator. The unbounded operators $K^{-2} \, , \, K^2 : \C A \to \C H_h$ play a special role in this paper. It can be proved that they are closable, with positive and selfadjoint closure. We will use the notation $\De_L = \overline{K^2}$, $\De_L^{-1} := \overline{K^{-2}}$ for the closures and refer to $\De_L$ as the \emph{left modular operator}. The left modular operator and the left modular action are related by the modular identity
\[
\si_L(x) = \De_L x \De_L^{-1} 
\]
for all $x \in \C A$.

The coordinate algebra $\C A$ has a vector space basis $\{t^l_{mn} \in \C A\,|\, 2l \in \nn \cup \{0\} \, , \, n,m = -l,\ldots,l\}$ consisting of matrix elements of its irreducible corepresentations, see \cite[Theorem 13, Section 4.2.5]{KS}. This basis for $\C A$ is actually an orthogonal basis for the surrounding Hilbert space $\C H_h$. We will use the notation $\{\xi^l_{mn}\}$ for the associated orthonormal basis for $\C H_h$ and refer to this basis as the \emph{corepresentation basis}. The unbounded operators $E,F,K$ and $K^{-1}$ associated with the generators of $\C U$ can be neatly described in terms of the corepresentation basis. Indeed we have that
\begin{equation}\label{eq:expefk}
\begin{split}
& E : \xi^l_{mn} \mapsto \sqrt{[l+n]_q[l-n+1]_q} \, \xi^l_{m,n-1} \\
& F : \xi^l_{mn} \mapsto \sqrt{[l - n]_q[l+n+1]_q}\, \xi^l_{m,n+1} \\
& K : \xi^l_{mn} \mapsto q^{-n} \xi^l_{mn}.
\end{split}
\end{equation}
Here we use the notation $[z]_q = \frac{q^z - q^{-z}}{q - q^{-1}}$ for the \emph{$q$-integer} associated with some complex number $z \in \cc$.

\section{The $q$-Dirac operator}\label{qdirac}
In this section we present the main construction of this paper and prove some basic properties. To be more precise, we obtain a triple $(\C A,\C H,D_q)$ consisting of the coordinate algebra for $\SU$, a Hilbert space $\C H$, and a selfadjoint unbounded operator $D_q$. As in the case of a spectral triple, both the $C^*$-algebra $\SU$ and the unbounded operator $D_q$ act on the Hilbert space. However, instead of the usual properties of a spectral triple, we only get bounded twisted commutators and furthermore, the resolvent lies in the compacts associated with a semifinite von Neumann algebra which does \emph{not} contain $\C A$. For convenience we will refer to our triple $(\C A, \C H, D_q)$ as a \emph{twisted spectral triple} even though we do not make a precise definition of this concept.
%

Let $\C H := \C H_h \op \C H_h \cong \C H_h \ot \cc^2$ and let $\pi := \pi \ot 1 : \SU \to \C L(\C H)$ denote the representation of $\SU$ given by diagonal left multiplication operators. 

We define the unbounded operator
\[
\De := \mat{cc}{K q^{-1/2} & 0 \\ 0 & K q^{1/2}} :
\C A \op \C A \to \C H
\]
and note that the closure $\T{cl}(\De)$ is selfadjoint and positive. For each $n \in \frac{1}{2}\zz$ we let $\C H_n \su \C H$ denote the eigenspace for $\T{cl}(\De)$ with eigenvalue $q^n$, and we remark that we have a Hilbert space decomposition $\C H \cong \T{cl}(\sum_{n \in \frac{1}{2}\zz} \C H_n)$. We will use the notation $\C M \su \C L(\C H)$ for the von Neumann algebra defined by $\C M := \{T \in \C L(\C H) \, | \, T(\C H_n) \su \C H_n \, \forall n \in \frac{1}{2}\zz\}$. The von Neumann algebra $\C M$ is semifinite and we fix the faithful normal semifinite trace $\Psi : \C M_+ \to [0,\infty]$ defined by
\[
\Psi(T) = \sum_{n \in \frac{1}{2}\zz}q^n\T{Tr}(T|_{\C H_n}) \, , \q T \in \C M_+.
\]
Here $\T{Tr}$ denotes the operator trace. We will use the notation $\C K(\C M,\Psi)$ for the associated $C^*$-algebra of compact operators. This is the smallest norm-closed ideal in $\C M$ which contains all projections $P$ with $\Psi(P) < \infty$.

We define the unbounded symmetric operator
\[
D := \mat{cc}{
(q^{-1}K^2 - 1)/(q - q^{-1}) & EK q^{1/2} \\
FK q^{-1/2} & (1 - q K^2)/(q - q^{-1})
} : \C A \op \C A \to \C H.
\]
We will apply the notation $D_q := \T{cl}(D)$ for the closure of $D$ and refer to this unbounded operator as the \emph{$q$-Dirac operator}.

We are now going to investigate the properties of the triple $(\C A,\C H,D_q)$, starting with selfadjointness and affiliation.

\begin{lemma}
The unbounded symmetric operator $D$ is essentially selfadjoint and the commutator $[D,\De] = 0$ is trivial. In particular, the $q$-Dirac operator $D_q$ is selfadjoint and affiliated with the von Neumann algebra $\C M$.
\end{lemma}
\begin{proof}
For each $l \in \frac{1}{2} \nn \cup \{0\}$ we let $\C A_l \su \C H_h$ denote the finite dimensional subspace with orthonormal basis $\{\xi^l_{mn}\, | \, m,n \in \{-l,\ldots,l\}\}$. We then note that the subspace $\C A_l \op \C A_l \su \C H_h$ is invariant for the symmetric unbounded operator $D$. In particular we can decompose $D$ as an infinite direct sum of selfadjoint operators. This proves that $D$ is essentially selfadjoint.

To prove that $[D,\De] = 0$ we compute as follows
\begin{equation}\label{eq:commvan}
D \De - \De D =
\mat{cc}{
0 & EK^2q - KEK \\ FK^2q^{-1} - KFK & 0
} = 0.
\end{equation}
Here we have applied the algebraic rules in $\C U$, see \eqref{eq:relU}.

To see that $D_q$ is affiliated with $\C M$ we note that each operator $T \in \C M'$ in the commutant can be obtained as a strong limit $T := \sum_{n \in S} \la_n P_n$ for some subset $S \su \frac{1}{2}\zz$. Here $P_n : \C H \to \C H$ denotes the projection with image $P_n(\C H) = \C H_n$ and $\{\la_n\}$ is a bounded sequence of complex numbers. Using the identity in \eqref{eq:commvan} it is not hard to verify that $T$ preserves the domain of $D_q$ and that the commutator $[D_q,T] : \T{Dom}(D_q) \to \C H$ is trivial.
\end{proof}

Our next concern is to show that the twisted commutators with $D_q$ and algebra elements result in bounded operators. The twist is given by the left modular action $\si_L \in \T{Aut}(\C A)$.

We will use the short notation $\pa_1 := \pa_{fk}, \, \pa_2:= \pa_{k^2 - 1}, \, \pa_3 := \pa_{ek} \in \T{End}(\C A)$. 
%

\begin{lemma}\label{l:twistcomm}
The twisted commutator $[D_q, x]_{\si_{L}}$, with twist given by the left action $\si_{L}$, extends to a bounded operator on $\C H$ for all $x \in \C A$. The twisted commutator can be presented explicitly as
\enveq{
\label{eq:comm}
D_q x - \si_L(x) D_q 
= \bma (q-q^{-1})^{-1} \pa_2(x) & q^{1/2}\pa_3(x) 
\\ q^{-1/2}\pa_1(x) & -(q-q^{-1})^{-1}\pa_2(x) \ema 
\in \C L(\C H)
}
for all $x \in \C A$.
\end{lemma}
\begin{proof}
Let $x \in \C A$ and let $\xi \in \C A \su \C H_h$. We then have the identities
\[
\begin{split}
& (\la K^2 -1)(x \cd \xi) 
- \si_L(x) \cd (\la K^2 -1)(\xi) \\
& \q = \si_L(x) \cd \xi - x \cd \xi
= \pa_2(x) \cd \xi
\end{split}
\]
for all $\la \in \cc$. Here we have used that $\pa_{k^2} = \si_L$ is an automorphism. Furthermore, using the derivational properties in \eqref{eq:twistdev} we obtain the identities
\[
\begin{split}
(EK)(x \cd \xi) - \si_L(x)\cd (EK)(\xi)
& = \pa_{ek}(x) \cd \xi  =  \pa_3(x) \cd \xi \\
(FK)(x \cd \xi) - \si_L(x)\cd (FK)(\xi)
& = \pa_{fk}(x) \cd \xi = \pa_1(x) \cd \xi.
\end{split}
\]
These computations prove that the identity in \eqref{eq:comm} is valid on the core $\C A \op \C A$ for the $q$-Dirac operator $D_q$, but this implies the statement of the lemma.

\end{proof}

We are now going to study the growth of the resolvent $(i + D_q)^{-1}$. We will use the trace ideal $\C L^1(\C M,\Psi)$ associated with the faithful normal semifinite trace $\Psi : \C M_+ \to [0,\infty]$. This is the ideal in $\C M$ defined by 
\[
\C L^1(\C M, \Psi) := \{x \in \C M \, | \, \Psi(|x|) < \infty\}.
\]
Equivalently, it consists of the elements in $\C M$ with Lebesgue integrable singular values (see \cite[Proposition $1.11$]{fack}). The trace ideal is contained in the compacts $\C K(\C M,\Psi)$ (see also \cite{fackkos}).

\begin{lemma}
The square $D^2 = \cas \De^2 : \C A \op \C A \to \C H$ recovers the quantum Casimir up to multiplication with the unbounded operator $\De^2 : \C A \op \C A \to \C H$. Furthermore, the resolvent $(i + D_q)^{-1} \in \C L^1(\C M,\Psi)$ is of trace class relative to $\Psi : \C M_+ \to [0,\infty]$.
\end{lemma}
\begin{proof}
In order to prove the relation $D^2 = \cas \De^2$ we start by noting that the diagonal and the anti-diagonal of $D$ anticommute, thus
\[
\begin{split}
& \mat{cc}{\frac{q^{-1}K^2 - 1}{q - q^{-1}} & 0 \\
0 & \frac{1 - qK^2}{q - q^{-1}}} 
\cd \mat{cc}{0 & EKq^{1/2}\\
FKq^{-1/2} & 0} \\
& \q = - \mat{cc}{0 & EKq^{1/2} \\
FKq^{-1/2} & 0}
\cd \mat{cc}{\frac{q^{-1}K^2 - 1}{q - q^{-1}} & 0 \\
0 & \frac{1 - qK^2}{q - q^{-1}}}.
\end{split}
\]
This observation allows us to compute as follows
\[
\begin{split}
D^2 & = 
\mat{cc}{\frac{(q^{-1}K^2 - 1)^2}{(q - q^{-1})^2} + EKFK & 0 \\
0 & \frac{(1 - qK^2)^2}{(q - q^{-1})^2} + FKEK
} \\
& = \mat{cc}{\frac{(q^{-1/2}K - q^{1/2} K^{-1})^2}
{(q - q^{-1})^{-2}} + EF & 0 \\
0 & \frac{(q^{-1/2}K^{-1} - q^{1/2}K)^2}{(q -q^{-1})^2}
+ FE } \\
& \q \cd \mat{cc}{q^{-1} K^2 & 0 \\ 0 & q K^2} \\
& = \cas \De^2.
\end{split}
\]
For the last identity we refer to the definition of the quantum Casimir in \eqref{eq:casimir}.

We now focus on proving that the resolvent of $D_q$ lies in the trace ideal. To this end we remark that the operator $D^2 = \cas \De^2$ can be decomposed as the strongly convergent sum
\[
D^2 = \sum_{l = 0,1/2,\ldots}\sum_{n=-l-1/2}^{l+1/2} q^{2n} [l+1/2]_q^2  P_{n,l} : \C A \op \C A \to \C H
\]
where $P_{n,l} : \C H \to \C H$ denotes the projection onto the finite dimensional subspace $\C H_n \cap (\C A_l \op \C A_l) \su \C H$. We recall that $\C H_n \su \C H$ denotes the eigenspace for $\T{cl}(\De)$ with eigenvalue $q^n$ whereas $\C A_l \su \C H_h$ is defined by the orthonormal basis $\{\xi^l_{mn}\, | \,n,m \in \{-l,\ldots,l\}\}$. In particular we have the expression
\[
(1 + D_q^2)^{-1/2} = 
\sum_{l = 0,1/2,\ldots}\sum_{n=-l-1/2}^{l+1/2} (1 + q^{2n} [l+1/2]_q^2 )^{-1/2} P_{n,l}
\]
for the absolute value of the resolvent as the limit of a bounded increasing sequence of positive operators.

To continue we remark that $P_{n,l} \in \C M$ and that $\Psi(P_{n,l}) \leq 2 (2l+1) q^n$. The result of the lemma now follows since the value
\[
\begin{split}
\Psi(1 + D_q^2)^{-1/2}
& \leq
\sum_{l = 0,1/2,\ldots}\sum_{n=-l-1/2}^{l+1/2}
2(2l+1)(1 + [l+1/2]_q^2 q^{2n})^{-1/2} q^n \\
& \leq
\sum_{l = 0,1/2,\ldots}2(2l+1)(2l+2) [l+1/2]_q^{-1} 
< \infty
\end{split}
\]
is finite.
\end{proof}

\begin{remark}
It is tempting to combine the notion of $\si$-spectral triples and semifinite spectral triples. That is, replacing the bounded operators in Definition \ref{sispectral} with a semifinite von Neumann algebra and using the compacts there to measure the size of the resolvent. This idea has turned out to be fruitful in the case of ordinary spectral triples, leading to an extension of the Connes-Moscovici local index theorem, see \cite{carphilrensukI, carphilrensukII,conmos}. We would like to emphasize that the resulting notion of a semifinite $\si$-spectral triple does \emph{not} cover the example of the present text. Indeed, our coordinate algebra $\C A$ is not contained in the von Neumann algebra $\C M$ and this is a serious obstacle for building a $KK$-theory class for the pair $SU_q(2)$-$\C K(\C M,\Psi)$. We refer to the recent paper \cite{RS} and the PhD-thesis of the second author \cite{senior} for a treatment of similar problems.
\end{remark}

We end this section by a study of Lipschitz-regularity in the sense of Definition \ref{sispectral}.

\begin{lemma}
The twisted commutator $[|D_q|,x]_{\si_L}$ extends to a bounded operator on $\C H$ for each $x \in \C A$.
\end{lemma}
\begin{proof}
We start by noting that we can restrict our attention to the first component in $\C H = \C H_h \op \C H_h$. Furthermore, since our twisted commutators are well-behaved with respect to products and adjoints we can assume that $x = a$ or $x = c$. We note that the unbounded operator $|D_q|$ acts on the first component of $\C H$ as $|D_q| : \xi^l_{mn} \mapsto [l+1/2]_q q^{-n-1/2} \xi^l_{mn}$. The generators $a$ and $c$ act as 
\[
\begin{split}
& a : \xi^l_{mn} 
\mapsto \al^+_{lmn} \xi^{l+1/2}_{m-1/2,n-1/2} 
+ \al^-_{lmn} \xi^{l-1/2}_{m-1/2,n-1/2} \\
& c : \xi^l_{mn} 
\mapsto \ga^+_{lmn} \xi^{l+1/2}_{m+1/2,n-1/2} 
+ \ga^-_{lmn} \xi^{l-1/2}_{m+1/2,n-1/2}
\end{split}
\]
where the coefficients satisfy the growth-constraints $|\al^+_{lmn}| \, , \, |\ga^+_{lmn}| \leq C_1 \cd q^{n+l}$ and $|\al^-_{lmn}| \, , \, |\ga^-_{lmn}| \leq C_2$ for all $l \in \{0,1/2,\ldots\}$ and $m,n \in \{-l,\ldots,l\}$. Here the constants $C_1,C_2 > 0$ are independent of the indices.

Using the above observations we can obtain the estimate
\[
\begin{split}
& \big|
\langle [|D_q|,x\,]_{\si_L}(\xi^l_{mn}),\xi^{l+1/2}_{m \pm 1/2,n-1/2} \rangle
\big| \\
& \q \leq C_1q^{n+l -1/2}\big| [l+1]_q q^{- (n-1/2)} - [l+1/2]_q q^{-n+1} \big| \\
& \q = C_1 q^{l-1/2}\big| [l+1]_q q^{1/2} - [l+1/2]_q q \big|
\end{split}
\]
which shows that the supremum 
\[
\T{sup}_{l,m,n} \big|\langle[|D_q|,x]_{\si_L}(\xi^l_{mn})
,\xi^{l+1/2}_{m \pm 1/2,n-1/2} \rangle\big| < \infty
\]
is finite. On the other hand we have the estimate
\[
\begin{split}
& \big|
\langle [|D_q|,x]_{\si_L}(\xi^l_{mn}),\xi^{l-1/2}_{m \pm 1/2,n-1/2}\rangle
\big| \\
& \q \leq C_2 q^{-1/2}\big| [l]_q q^{-(n-1/2)} - [l+1/2]_q q^{-n +1} \big| \\
& \q = C_2 q^{-n}\big| (q^l - q^{-l} - q^{l+1} + q^{-l})/(q - q^{-1}) \big| \\
& \q = C_2 q^{l-n} \cd (q-1)/(q - q^{-1})
\end{split}
\]
which shows that the supremum
\[
\T{sup}_{l,m,n}
\big|\langle[|D_q|,x]_{\si_L}(\xi^l_{mn}),
\xi^{l-1/2}_{m \pm 1/2,n-1/2}\rangle\big| < \infty
\]
is finite. This ends the proof of the lemma.
\end{proof}

\section{Spectral dimension and the $q$-zeta function}\label{s:qzet}
In this section we will investigate the $q$-zeta function
\[
z \mapsto \Psi(|D_q|^{-z}) := \ze_q(z)
\]
associated with our twisted spectral triple $(\C A,\C H,D_q)$ and the semifinite trace $\Psi : \C M_+ \to [0,\infty]$. We shall then see that our $q$-zeta function is holomorphic for all $z \in \cc$ with $\T{Re}(z) > 1/2$ and that it extends to a meromorphic function on the whole complex plane with poles of order two at $1/2 - \nn \cup \{0\} + i \cd 2\pi \log(q)^{-1} \cd \zz$. We remark that the $q$-zeta function which we obtain is related to the zeta functions studied in \cite{kawwakyam,uennis} among others. In particular our proof of the existence of the meromorphic extension follows a well known argument. We also use the opportunity to refer to the paper \cite{krawag}, where a similar zeta function is studied for the standard quantum sphere.

\begin{prop}\label{merext}
The $q$-zeta function $z \mapsto \Psi(|D_q|^{-z}) := \ze_q(z)$ is holomorphic for all $z \in \cc$ with $\E{Re}(z) > 1/2$ and has a meromorphic extension to the whole complex plane with poles of order two at $1/2 - \nn \cup \{0\} + i \cd 2\pi \log(q)^{-1} \cd \zz$. The meromorphic extension is given by the explicit formula
\[
\ze_q(z) =
(q^{-1} - q)^z
\cd \sum_{j=0}^\infty {z + j - 1 \choose j} \cd 
\frac{(q^{z-1/2+j} + q^{j+1/2}) \cd (1 - q^{2j+z})}
{(1-q^{z-1/2 + j})^2(1 - q^{j+1/2})^2}
\]
for all $z \in \cc$.
\end{prop}
\begin{proof}
Let us fix some $z \in \cc$ with $\T{Re}(z) > 1/2$. We then have that
\[
\begin{split}
\Psi(|D_q|^{-z})
& = \sum_{l=0,1/2,\ldots} \sum_{n=-l-1/2}^{l+1/2} q^{-nz}[l+1/2]_q^{-z} \Psi(P_{n,l}) \\
& = (q^{(z-1)/2} + q^{(1-z)/2}) \cd
\sum_{l=0,1/2,\ldots} \sum_{n=-l}^l (2l+1) q^{-nz +n}[l+1/2]_q^{-z}.
\end{split}
\]
To continue from this point we note that $\sum_{n=-l}^l q^{n(1-z)} = [2l+1]_{q^{(1-z)/2}}$ where by definition
\[
[k]_w := \fork{ccc}{
(w^{-k} - w^k)/(w^{-1} - w) & & w \neq 0 \\
k & & w = 0
}
\]
for all $k \in \nn$ and all $w \in \cc$. Furthermore, we have that
\[
\begin{split}
[l+1/2]_q^{-z}
& = (q^{-1} - q)^z \cd q^{(l+1/2)z}(1 - q^{2l+1})^{-z} \\
& = (q^{-1} - q)^z \cd q^{(l+1/2)z}
\sum_{j=0}^\infty {z+j-1 \choose j} q^{(2l+1)j}
\end{split}
\]
where the series converges absolutely. Now, let us fix some $j \in \nn \cup \{0\}$. It can then be proved that
\[
\sum_{k=1}^\infty k q^{k(j + z/2)}[k]_{q^{(1-z)/2}}
= \frac{q^{z/2+j} \cd (1 - q^{2j+z})}
{(1-q^{z-1/2 + j})^2(1 - q^{j+1/2})^2}
\]
where the series converge absolutely. Indeed, this follows from basic algebraic manipulations together with the summation formulae
\[
\begin{split}
\sum_{k=1}^\infty k q^{w \cd k} = \frac{q^w}{(1 - q^w)^2} \qq
\sum_{k=1}^\infty k^2 q^{w \cd k} = \frac{q^w(1 + q^w)}{(1 - q^w)^3}
\end{split}
\]
which are valid for any complex number $w$ with $\T{Re}(w) > 0$.

Combining the above computations we obtain the expression
\[
\Psi(|D_q|^{-z}) 
= (q^{-1} - q)^z
\cd \sum_{j=0}^\infty {z + j - 1 \choose j} \cd 
\frac{(q^{z-1/2+j} + q^{j+1/2}) \cd (1 - q^{2j+z})}
{(1-q^{z-1/2 + j})^2(1 - q^{j+1/2})^2}
\]
for the $q$-zeta function. But this proves the claim of the proposition.
%
\end{proof}

\begin{remark}
It follows from Proposition \ref{merext} that our twisted spectral triple has spectral dimension $1/2$, in the sense that the bounded operator $(1 + D_q^2)^{-r/2}$ lies in the trace ideal $\C L^1(\C M,\Psi)$ for all $r > 1/2$ but not for $r =1/2$, see \cite[Definition $4.1$]{carphilrensukI}. We would like to remark that this particular spectral dimension is linked with the choice of the trace $\Psi : \C M_+ \to [0,\infty]$. To be more precise, replacing $\Psi$ with the trace $\Psi_s : \C M_+ \to [0,\infty]$ defined by
\[
\Psi_s(T) = \sum_{n \in \frac{1}{2}\zz} q^{sn}\E{Tr}(T|_{\C H_n}) \q \T{for some } s > 0
\]
we get that the corresponding spectral dimension becomes $s/2$. Or in other words, by changing the trace we can obtain any strictly positive spectral dimension.
\end{remark}





\section{The triviality of local twisted Hochschild $3$-cocycles}
\label{sec_ch2_trivial}

In this section we will show that any "local" twisted Hochschild $3$-cocycle constructed using twisted commutators with $D_q$ is trivial, where $D_q$ is the $q$-Dirac operator from the twisted spectral triple $(\C A, \C H, D_q)$ defined in Section \ref{qdirac}. Our local twisted Hochschild cocycles are analogous to the Hochschild class of the Chern character from the local index formula, however we replace commutators with $D_q$ by twisted commutators with $D_q$. We refer to the formula for the local Hochschild cocycle given in \cite[Proposition $3.10$]{conmosc}.
%

The case of twisted Hochschild $3$-cocycles is of particular interest in view of the computation of twisted Hochschild homology carried out by T. Hadfield and U. Kr\"ahmer in \cite{hadkrah}, and because the classical Lie group $SU(2)$ has dimension $3$. The computations of Hadfield and Kr\"ahmer reveal a fundamental class in the third twisted Hochschild homology group associated with the modular automorphism of the Haar-state. Our results show that the twisted spectral triple does \emph{not} detect this class by means of a local Hochschild cocycle construction.

Let us recall the definition of twisted Hochschild cohomology, see \cite{kustmurtu}. Let $\C B$ be a unital algebra over $\cc$ and let $\si \in \T{Aut}(\C B)$ be some automorphism of $\C B$. For each $n \in \nn \cup \{0\}$ we let $C^n(\C B,\si)$ denote the vector space defined by
\[
C^n(\C B, \si) := \big\{ \varphi \in \T{Hom}_{\cc}(\C B^{\ot (n+1)},\cc) \, | \, 
\varphi \ci \si = \varphi
\big\}.
\]
Here $\si : \C B^{\ot (n+1)} \to \C B^{\ot (n+1)}$ denotes the isomorphism defined by $\si(b_0\olo b_n) = \si(b_0) \olo \si(b_n)$. We refer to $C^n(\C B,\si)$ as the \emph{twisted Hochschild $n$-cochains}. The twisted Hochschild cochains can be equipped with a coboundary operator $b_\si : C^n(\C B,\si) \to C^{n+1}(\C B,\si)$ defined by
\[
\begin{split}
b_{\si}(\varphi)(b_0\olo b_{n+1}) & =
(-1)^{n+1} \varphi(\si(b_{n+1})b_0 \olo b_n) \\
& \q + \sum_{i=0}^n (-1)^i \varphi(b_0 \olo b_i b_{i+1} \olo b_{n+1}).
\end{split}
\]
We define the \emph{twisted Hochschild cohomology} as the cohomology of the complex
\[
\begin{CD}
0 @>>> C^0(\C B,\si) @>b_\si>> C^1(\C B,\si) @>b_\si>> \ldots.
\end{CD}
\]
The cohomology group in degree $n$ will be denoted by $HH^n(\C B,\si)$. We will refer to a twisted Hochschild cochain $\varphi \in C^n(\C B,\si)$ with $b_\si(\varphi) = 0$ as a \emph{twisted Hochschild cocycle}. A twisted Hochschild cochain which lies in the image of $b_\si$ is called a \emph{twisted Hochschild coboundary}. Clearly each twisted Hochschild cocycle $\varphi$ defines a cohomology class $[\varphi] \in HH^n(\C B,\si)$ and this class is trivial if and only if $\varphi$ is a twisted Hochschild coboundary.
%

Let us describe a general method for constructing twisted Hochschild cocycles. Suppose that $\tau : \C B \to \cc$ is a linear functional which is a twisted trace for some automorphism $\al \in \T{Aut}(\C B)$. Thus we have the identity $\tau(x\cdot y) = \tau(\al(x) \cdot y)$ for all $x,y \in \C B$. Furthermore, suppose that $d_1,\ldots,d_n \in \T{End}(\C B)$ are twisted derivations on $\C B$ with respect to some automorphism $\te \in \T{Aut}(\C B)$. Thus we have the identity $d_{i}(x y) = d_{i}(x) y + \te(x) d_{i}(y)$ for all $i = 1,\ldots,n$ and all $x,y \in \C B$. We define the automorphism $\si := \al \ci \te^{-n} \in \T{Aut}(\C B)$.

On top of these conditions we will assume that the automorphisms commute, thus $\te \ci \al = \al \ci \te$. And furthermore, we assume that the twisted trace and the derivations satisfy the invariance properties $\tau(\si(x)) = \tau(x)$ and $d_1(\si(x_1))\clc d_n(\si(x_n)) =\si( d_1(x_1) \clc d_n(x_n))$ for all $x,x_1,\ldots,x_n \in \C B$.

We can then form a twisted Hochschild $n$-cochain $\varphi \in C^n(\C B,\si)$ using the formula
\[
\varphi : (x_0,\ldots,x_n) \mapsto
\tau\big( x_0 d_1(\te^{-1}(x_1)) \clc
d_n( \te^{-n}(x_n)) \big)
\]
for $x_i \in \C B$. Here the twist is given by $\si =\al \ci \te^{-n} \in \T{Aut}(\C B)$. Remark that $\varphi \ci \si = \varphi$ by the commutativity of the automorphisms and the invariance conditions on the twisted trace and the derivations. Furthermore, using that the endomorphisms $d_i \in \T{End}(\C B)$ are twisted derivations and that $\tau : \C B \to \cc$ is a twisted trace we get that the twisted Hochschild coboundary $b_\si(\varphi) = 0$ is trivial. See for example the computation in \cite[Proposition $3.10$]{conmosc}. In particular we get a class $[\varphi] \in HH^n(\C B,\si)$ in the $n^{\T{th}}$ twisted Hochschild cohomology group.

The following lemma shows that $\varphi$ is a twisted Hochschild coboundary when one of the derivations is a twisted inner derivation.

\begin{lemma}\label{l:gentriv}
Suppose that $d_j \in \E{End}(\C B)$ agrees with the inner twisted derivation defined by some $Y \in \C B$, thus $d_j(x) = Y x - \te(x) Y$ for all $x \in \C B$. Suppose furthermore that the derivations satisfy the invariance property
\begin{equation}\label{eq:invarder}
\begin{split}
& \si\big(d_1(x_1)\clc d_{j-1}(x_{j-1})
\cd Y \cd d_{j+1}(x_j)\clc d_n(x_{n-1})\big) \\
& \q = d_1(\si(x_1)) \clc d_{j-1}(\si(x_{j-1})) \cd Y \cd
d_{j+1}(\si(x_j))\clc d_n(\si(x_{n-1}))
\end{split}
\end{equation}
for all elements $x_1,\ldots,x_{n-1} \in \C B$. The twisted Hochschild cocycle $\varphi$ is then a twisted Hochschild coboundary. In particular the class $[\varphi] \in HH^n_\si(\C B)$ is trivial.
\end{lemma}
\begin{proof}
 We define the twisted Hochschild $(n-1)$-cochain $\varphi_j : \C B^{\ot (n-1)} \to \cc$ by the formula
\[
\begin{split}
\varphi_j : (x_0,\ldots,x_{n-1}) & \mapsto 
\tau \big(x_0 d_1(\te^{-1}(x_1)) \clc 
d_{j-1}(\te^{-(j - 1)}(x_{j-1})) \\
& \qq \cd Y d_{j+1}(\te^{-(j+1)}(x_j)) 
\clc d_n(\te^{-n}(x_{n-1})) \big).
\end{split}
\]
We remark that the invariance property $\varphi_j \ci \si = \varphi_j$ follows from the condition in \eqref{eq:invarder} together with the commutativity of the automorphisms and the invariance condition on the twisted trace. Furthermore, using that the endomorphisms $d_1,\ldots,d_{j-1},d_{j+1},\ldots,d_n \in \T{End}(\C B)$ are twisted derivations for $\te \in \T{Aut}(\C B)$ and that $\tau$ is a twisted trace with respect to $\al$, it is not hard to see that $b_\si(\varphi_j) = (-1)^j \varphi$. But this proves the statement of the lemma.
\end{proof}


We now specialise to the case of the coordinate algebra $\C A$ for $\SU$. Let $\al \in \T{Aut}(\C A)$ denote the automorphism defined by $\al : \xi^l_{mn} \mapsto \la^m \mu^n \xi^l_{mn}$ for some fixed $\la, \mu \in \cc\setminus \{0\}$. We suppose that $\tau : \C A \to \cc$ is a twisted trace for the automorphism $\al$.

We extend the functional $\tau$ to the matrix algebra $M_2(\C A) = \C A \ot M_2(\cc)$ by putting $\tau = \tau \ot \T{Tr}$ where $\T{Tr} : M_2(\cc) \to \cc$ denotes the usual trace on two by two matrices. We then define the functional $\varphi : \C A^{\ot 4} \to \cc$ by the formula
\begin{equation*}
\varphi(x_0,\ldots,x_3) = \tau(x_{0} [D_q, \si_{L}^{-1}(x_{1})]_{\si_{L}} [D_q, \si_{L}^{-2}(x_{2})]_{\si_{L}} [D_q, \si_{L}^{-3}(x_{3})]_{\si_{L}}).
\end{equation*}
We remark that the above functional looks like the local Hochschild cocycle defined by Connes and Moscovici in \cite[Proposition $3.10$]{conmosc} except that their twisted Dixmier trace has been replaced by a general twisted trace $\tau : M_2(\C A) \to \cc$. Note that $D_q$ is the $q$-Dirac operator which we constructed in Section \ref{qdirac}.

Using the computations of twisted commutators given in Lemma \ref{l:twistcomm} we obtain the expression
\[
\varphi = \sum_{s \in \Si_3} \frac{\T{sgn}(s)}{ q - q^{-1}} \varphi_s
\]
where $\Sigma_3$ is the permutation group on three letters and $\varphi_s$ is defined by
\[
\varphi_s : (x_0,\ldots,x_3) \mapsto 
\tau\big( x_0 \pa_{s(1)}(\si_L^{-1}(x_1))
\pa_{s(2)}(\si_L^{-2}(x_2)) 
\pa_{s(3)}(\si_L^{-3}(x_3)) \big).
\]

The next result is a consequence of our triviality result for twisted Hochschild cocycles associated with twisted traces and inner twisted derivations. It implies that $\varphi$ defines a trivial class in the third twisted Hochschild cohomology group $HH^3(\C A,\si)$ where $\si := \al \ci \si_L^{-3} \in \T{Aut}(\C A)$.
%
%

\begin{prop}
Suppose that the $\al$-twisted trace $\tau : \C A \to \cc$ is $\si_L$-invariant, thus $\tau(\si_L(x)) = \tau(x)$ for all $x \in \C A$. Then the functional $\varphi_s : \C A^{\ot 4} \to \cc$ is a twisted Hochschild coboundary for the automorphism $\si = \al \ci \si_L^{-3} \in \E{Aut}(\C A)$ for all $s \in \Si_3$. In particular we get that the functional $\varphi$ is a twisted Hochschild coboundary.
\end{prop}
\begin{proof}
We remark that $\pa_2 \in \T{End}(\C A)$ is a twisted inner derivation with respect to $\si_L$. Indeed, we have that 
\[
\pa_2(x) = \si_L(x) - x = [-1,x]_{\si_L} 
\]
for all $x \in \C A$. Furthermore, it follows from \eqref{eq:twistdev} that $\pa_1 = \pa_{fk} \T{ and } \pa_3 = \pa_{ek} \in \T{End}(\C A)$ are twisted derivations for the left modular automorphism $\si_L$.

Next, we note that $\al \ci \si_L = \si_L \ci \al$ and $\tau(x) = \tau(\si(x))$ for all $x \in \C A$. Furthermore, we have the identities $\si(\pa_1(x)) = \mu q^6 \pa_1(\si(x))$ and $\si(\pa_3(y)) = \mu^{-1} q^{-6}\pa_3(\si(y))$ for all $x,y \in \C A$. See \eqref{eq:expefk}.

The result of the proposition now follows by an application of Lemma \ref{l:gentriv}.
\end{proof}

\begin{remark}
Using similar methods it can be proved that the formula
\[
(x_0,\ldots,x_{2k+1}) \mapsto \tau(x_0[D_q,\si_L^{-1}(x_1)]_{\si_L}\clc 
[D_q,\si_L^{-2k-1}(x_{2k+1})]) 
\]
defines a twisted Hochschild $(2k+1)$-coboundary for all $k \in \nn \cup \{0\}$. This leaves the problem of understanding the homological dimension of the twisted spectral triple $(\C A, \C H, D_q)$ open.
\end{remark}

\end{document}